\documentclass{adamjoucc}
\usepackage{lastpage}
\usepackage{color}
\usepackage{amssymb}
\usepackage{mathtools}
\usepackage{tikz}
\usepackage{chemfig}
\usepackage{ogonek} % for polish names
\usetikzlibrary{calc, positioning, quotes}

\makeatletter
\def\moverlay{\mathpalette\mov@rlay}
\def\mov@rlay#1#2{\leavevmode\vtop{%
    \baselineskip\z@skip \lineskiplimit-\maxdimen
    \ialign{\hfil$\m@th#1##$\hfil\cr#2\crcr}}}
\newcommand{\charfusion}[3][\mathord]{
  #1{\ifx#1\mathop\vphantom{#2}\fi
    \mathpalette\mov@rlay{#2\cr#3}
  }
  \ifx#1\mathop\expandafter\displaylimits\fi}
\DeclareRobustCommand\bigop[1]{%
  \mathop{\vphantom{\sum}\mathpalette\bigop@{#1}}\slimits@
}
\newcommand{\bigop@}[2]{%
  \vcenter{%
    \sbox\z@{$#1\sum$}%
    \hbox{\resizebox{\ifx#1\displaystyle.9\fi\dimexpr\ht\z@+\dp\z@}{!}{$\m@th#2$
}}%
  }%
}
\makeatother
\DeclareMathOperator{\cut}{cut}

\journalyear{2025}
\journalpages{1}{\pageref{LastPage}}
\evenrunninghead{N.\ Domschke, T.\ Gatter, R. Golnik, P.F.\ Stadler}
\oddrunninghead{N.\ Domschke, T.\ Gatter,  R. Golnik, P.F.\ Stadler:
  Relevant Cuts}
\oddrunninghead{Relevant Cuts in Graphs}
\makeatletter
\def\@journalname{Preprint}%{Discrete Mathematical Chemistry}
\def\@issn{ }
\def\@journalnumber{}%{(for review)}

\makeatother

\definecolor{orange}{rgb}{1,0.5,0}
 \newcommand{\REV}[1]{\begingroup\color{blue}#1\endgroup}

\begin{document}

\begin{frontmatter}

\titledata{A Short Note on Relevant Cuts}{}

\authordata{Nico Domschke}{Bioinformatics Group, Department of Computer
  Science, Leipzig University, H{\"a}rtelstra{\ss}e 16–18, D-04107 Leipzig,
  Germany}{dnico@bioinf.uni-leipzig.de}{}

\authordata{Thomas Gatter}{Bioinformatics Group, Department of Computer
  Science, Leipzig University, H{\"a}rtelstra{\ss}e 16–18, D-04107 Leipzig,
  Germany}{thomas@bioinf.uni-leipzig.de}{}
    
\authordata{Richard Golnik}{Bioinformatics Group, Department of Computer
  Science, Leipzig University, H{\"a}rtelstra{\ss}e 16–18, D-04107 Leipzig,
  Germany}{richard@bioinf.uni-leipzig.de}{}    
        
\authordata{The Students of The Advanced Bioinformatics 
Course 2024}{Department of Computer
  Science, Leipzig University, Augustusplatz 10, D-04109 Leipzig,
  Germany}{}{}
  
\authordata{Peter F.\ Stadler}{Bioinformatics Group,
  Department of Computer Science \&
  Interdisciplinary Center for Bioinformatics \&
  Center for Scalable Data Analytics and Artificial Intelligence 
  Dresden/Leipzig \& 
  School of Embedded Composite Artificial Intelligence,
  Leipzig University, H{\"a}rtelstra{\ss}e 16–18,
  D-04107 Leipzig, Germany;
  Max Planck Institute for Mathematics in the Sciences,
  Inselstra{\ss}e 22, D-04103 Leipzig, Germany;
  Department of Theoretical Chemistry,
  University of Vienna, W{\"a}hringerstra{\ss}e 17,
  A-1090 Wien, Austria;
  Facultad de Ciencias, Universidad National de Colombia;
  Bogot{\'a}, Colombia;
  Santa Fe Institute, 1399 Hyde Park Rd., Santa Fe 
  NM 87501, USA}{studla@bioinf.uni-leipzig.de}{}

\keywords{
  Cut Space;
  Matroid;
  Gomory-Hu tree;
  Enumeration;
  %  Complex Reactions;
}
\msc{05C21,05C50} 

\begin{abstract}
  The set of relevant cuts in a graph is the union of all minimum weight
  bases of the cut space. We show that a cut is relevant if and only if it
  is a minimum weight cut between two distinct vertices. Moreover, we give
  a characterization in terms of Picard-Queyranne Directed Acyclic
    Graphs that can be used to accelerate the enumeration of the relevant
  cuts.  Finally, we perform an experimental evaluation by comparing with
state-of-the-art algorithms.
\end{abstract}

\end{frontmatter}

\section{Introduction}

Cuts in graphs are closely related to cycles. The cycle space and the cut
space form orthogonal subspaces of the $|E|$-dimensional vector space over
GF(2) spanning the edges of a graph $G = (V, E)$. Nevertheless, graphs can contain edge sets
called bi-cycles that are both cycles and cuts \cite{Rosenstiehl:78}. The
cycle space, and in particular bases of the cycle space comprising short
cycles, play an important role in cheminformatics
\cite{Downs:89,Kim:09,May:14} and structural biology
\cite{Leydold:98a,Lemieux:06}. In particular, minimum cycle bases composed
of cycles with minimal total length or edge weight, and their union, known
as the set of \emph{relevant cycles} \cite{Vismara:97}, are of use
\cite{Berger:04a}.

Like cycle bases, cut bases can also be defined in directed graphs, where
they are usually referred to as co-cycle bases \cite{Kavitha:12}. In this
case, co-cycles are considered as vectors over $\mathbb{Q}$. Co-cycle bases
have found applications in particular in crystallography
\cite{DelgadoFriedrichs:03,Xiao:23}. 

Minimum cut bases in undirected graphs are studied in detail in
\cite{Bunke:10}, their counterparts in directed graphs were analyzed in
\cite{Kavitha:08,Kavitha:12}, highlighting the special role of Gomory-Hu
trees. These are (not necessarily spanning) trees on $V(G)$ with edge
weights $\lambda_e$ such that for each pair $s,t\in V(G)$, the minimum
weight $\lambda_e$ along the path between $s$ and $t$ in $T$ equals the
weight of a minimum $s,t$-cut \cite{Gomory:61}. Since Gomory-Hu trees exist
for all graphs, certain sets of $n-1$ distinct minimal cuts already provide
a minimal cut for any pair of distinct vertices. As shown in
\cite{Bunke:10,Kavitha:08}, the cuts defined by a Gomory-Hu tree form
minimal cut (co-cycle) bases. Moreover, they are weakly fundamental, i.e.,
their cuts can be ordered in such a way that every cut contains an edge
that is not present in any of its predecessors.

The construction principle of Gomory-Hu trees ensures that the associated
cuts are nested, i.e., for any two distinct cuts $(A,\bar A)$ and $(B,\bar
B)$, one of the intersections $A\cap B$, $A\cap \bar B$, $\bar A\cap B$ and
$\bar A\cap \bar B$ is empty. A simplified algorithm \cite{Gusfield:90},
however, alleviates this restriction, suggesting that minimal cut bases are
more general than the Gomory-Hu bases. 

Beyond the interest in relevant cuts implied by the analogy to relevant
cycles, these are also of practical interest\REV{,} e.g.\REV{,}\ for the design of
well-behaved crossover operators in graph-based genetic
algorithms. For example, recently, cut-and-join crossover operators for chemical graphs
have been investigated and benchmarked in detail for search tasks in
``chemical space'', demonstrating their efficacy \cite{Domschke:25a}.

In the following sections we show that relevant cuts coincide
    with the cuts that are weight-minimal for some pair of vertices; these
    results follows from previously known properties of minimal cut
    bases (see Sections 2 and 3). Then, we use these insights to construct novel algorithms for
    computing a cut basis and the set of relevant cuts of arbitrary
  graphs (Section 4). Moreover, we benchmark pre-existing and novel algorithms on
  chemical and random graphs (Section 5). Finally, we draw our conclusions (Section 6).

\section{Preliminaries}

Throughout, $G=(V,E,w)$ is an undirected connected graph with edge weights
$w:E\to\mathbb{R}$. In the context of chemical applications,
integer weights corresponding to the order of a chemical bond will be most
relevant. A cut is a proper bipartition $(V',V'')$ of the vertex set
$V(G)$, i.e., $V'\ne \emptyset$, $V''\ne \emptyset$, $V'\cap
V''=\emptyset$, and $V'\cup V''=V$. A cutset $D\coloneqq
\cut(V',V'')=\{xy\in E| x\in V', y\in V''\}$.  The weight of a cutset is
$|D|\coloneqq\sum_{xy\in D} w(xy)$. In the following we will usually write
$\partial V'=\partial V''=\cut(V',V'')=\cut(V'',V')$. In particular, the
edges incident to each vertex $x\in V$ form a cutset $\partial\{x\}$. A
cutset $D$ that does not include another cutset $D'$ of
$G$ as a proper subset is called a \emph{bond} of $G$. Clearly, each
$\partial\{x\}$ is a bond. Let $T$ be a spanning tree of $G$, and $e$ an edge
of $T$. Removal of $e$ decomposes $T$ into exactly two connected components
with vertex sets $V_e'$ and $V_e''$. By construction the induced subgraphs
$G[V_e']$ and $G[V_e'']$ are connected.  The cut $\partial V_e'$ is also a
bond.

The cutsets of $G$ form a vector space $\mathfrak{D}\coloneqq
\mathfrak{D}(G)$ over GF(2) called the \emph{cut space}. We write $\oplus$
for the vector addition and note that this corresponds to the symmetric
difference of edge sets. Thus we have $D\oplus D=\varnothing$.  It is not
difficult to check, moreover, that
\begin{equation}
  \partial W = \bigoplus_{x\in W} \partial\{x\}
\end{equation}
for any $W \subseteq V$. For connected graphs, $\mathfrak{D}$ has dimension $n-1$, where $n=|V|$ is
the number of vertices in $G$. In particular, each of the sets
$\mathcal{D}_{x_0}\coloneqq \{\partial\{x\}|x\in V\setminus\{x_0\}\}$ for
any $x_0\in V$, and each of the sets $\mathcal{D}_T\coloneqq \{\partial
V_e'|e\in E(T)\}$ for any spanning tree $T$ of $G$ is basis of the cut
space, see e.g.\ \cite{Biggs:93}.

A cut $(V',\bar{V}')$ is an $s,t$-cut if $s\in V'$ and $t\in\bar{V}'$.
Clearly, the cut space $\mathfrak{D}$ contains an $s,t$-cut for every pair
$s\ne t$ of distinct vertices. Less obviously, this is also true for
  every basis $\mathcal{D}$ of $\mathfrak{D}$, see, e.g.,\ \cite{Bunke:10}.
  For completeness, we include a short proof of this statement, which will
  be used repeatedly throughout this contribution:
\begin{lemma}
  Let $G$ be a connected graph, $s,t\in V(G)$ distinct vertices, and  
  $\mathcal{D}$ a cut basis. Then $\mathcal{D}$ contains an $s,t$-cut.
  \label{lem:Dst} 
\end{lemma}
\begin{proof}
  Let $\mathcal{D}=\{ \partial U_i| i=1,\dots,n-1\}$ be a cut-basis
  comprising cuts $(U_i,V\setminus U_i)$ and assume that $\mathcal{D}$
    does not contain an $s,t$-cut. Without loosing generality, we may
    therefore assume that $s,t\in U_i$ for all $i$. Consider the set
    $Z\coloneqq\bigcap_i U_i$. We have $s,t\in Z$, and moreover, any
    cut that can be obtained as an $\oplus$-sum of cuts in $\mathcal{D}$ if
    of the form $(V',\bar{V'})$ with $Z\subseteq V'$ or
    $Z\subseteq\bar{V'}$ since none of the cuts $\partial U_i$ contains an
    edge with both endpoints in $Z$. Hence, note of the cuts $(V',\bar V')$
    is an $s,t$-cut, i.e., the $\oplus$-operation on $\mathcal{D}$ cannot
    generate any $s,t$-cut. The cut space, however, contains in particular
    the $s,t$-cut $\partial\{s\}$, contradicting our assumption that the
    $\mathcal{D}$ is a cut basis.
\end{proof}
  
\begin{definition}
  A basis $\mathcal{D}$ of the cut space $\mathfrak{D}$ is a \emph{minimum
  cut basis} if $\sum_{D\in\mathcal{D}} |D|$ is minimum.
\end{definition}

\begin{definition} 
Let $G=(V,E,w)$ be a weighted graph and let $T^*=(V,F,\lambda)$ be an
edge-weighted tree on the same vertex set. Then $T^*$ is a \emph{Gomory-Hu
tree} for $G=(V,E,w)$ if
\begin{itemize}
  \item[(G1)] for every $uv\coloneqq e\in F$, the bipartition of $V$
    induced by removing $e$ from $F$ defines a cut of $G$, which we denote
    by $(V_e,\bar V_e)$;
  \item[(G2)] for every $uv\coloneqq e\in F$, the weight $\lambda_e$ of $e$
    is the weight of a minimum $s,t$-cut;
  \item[(G3)] for all $s,t\in V$, the weight of a minimum $s,t$-cut in $G$
    equals the minimum weight $\lambda_e$ along the unique path $P_{st}$
    connecting $s$ and $t$ in $T^*$.
\end{itemize}
\end{definition} 
A key result in the theory of flow networks is
\begin{proposition} \cite{Gomory:61}
  Every weighted graph $G=(V,E,w)$ has a Gomory Hu tree
  $T^*=(V,F,\lambda)$.
\end{proposition}
\begin{proposition} \cite{Bunke:10}
  The cutsets $\partial{V_e}$ induced by the edges $e\in F$ of a Gomory-Hu
  tree $T^*$ for $G$ form a minimum cut basis of $G$.
\end{proposition}

While there may be exponentially many minimal $s,t$-cuts for any such fixed node pair,
Picard and Queyranne \cite{Picard:80} introduced a directed acyclic graph (as a shorthand, we write a PQ-DAG $Q_{s,t}$) that acts as a compact representation and may nevertheless be derived in polynomial time.
It is obtained from $G$ by contracting each connected component of the residual flow network to a single \emph{node} for an arbitrary maximum $s,t$-flow. 
Accordingly, by construction, vertices in $Q_{s,t}$ represent equivalence classes of vertices w.r.t. the strong connected relationship in the residual flow network, hence sets  of nodes from $G$.
The ``closed'' node sets in the DAG $Q_{s,t}$, i.e., those without out-edges, correspond to the minimal $s,t$-cuts.

As shown in \cite{Gusfield:93}, given PQ-DAG representations $Q_{u,v}$, corresponding to all minimal $u,v$-cuts, for the edges of the Gomory-Hu tree $uv \in T^*$, any PQ-DAG $Q_{s,t}$, where $st$ is not an edge in  $T^*$, can be computed as a recombination of PQ-DAGs corresponding to edges along the unique path $P_{st}$ connecting $s$ and $t$ in $T^*$, improving algorithmic complexity compared to direct construction.

\section{Relevant Cuts} 

The notion of \emph{relevant cycles} was studied for cycles bases of
undirected graphs in \cite{Vismara:97}.  Interestingly, the notion was
proposed already in 1971 by Morris Plotkin in the context of chemical ring
perception \cite{Plotkin:71}. More generally, the concept of \emph{relevant
elements} pertains to matroids with a weight function in general and reads
as follows for the special case of cut spaces:
\begin{definition}
  A cut $D\in\mathfrak{D}$ is \emph{relevant} if there is a minimum cut
  basis $\mathcal{D}$ of $\mathfrak{D}$ such that $D\in\mathcal{D}$.
\end{definition}
The following simple characterization of relevant elements holds for
matroids in general \cite{Vismara:97}:
\begin{lemma}
  \label{lem:Vismara}
  A cut $D\in\mathfrak{D}$ is relevant if and only if it cannot be written
  as an $\oplus$ sum of shorter cuts $D_1,\dots,D_k\in\mathfrak{D}$, i.e.,
  such that $|D_i|<|D|$ for $1\le i\le k$.
\end{lemma}
\begin{proof}
  Since $\mathfrak{D}$ is a vector space it is also matroid. Suppose $D$ is
  contained in the minimum cut basis $\mathcal{D}$.  If $D$ can be written
  as the $\oplus$-sum of shorter cuts $D_1$ through $D_k$, it could be
  exchanged with one of these cuts, contradicting minimality of
  $\mathcal{D}$. This is justified as follows: We know that at least one cut $D_1$ through $D_k$ is not element of $\mathfrak{D}$, otherwise $D\notin\mathfrak{D}$ by linear independence. If exactly one cut $D_i$ is not an element of $\mathfrak{D}$, $D_i$ and $D$ can be exchanged. If multiple cuts are not in $\mathfrak{D}$, we know they can be constructed as $\oplus$-sums of yet shorter cuts, and we apply this argument recursively.
  
  Now suppose that $D$ cannot be written as a sum of shorter
  cuts. Then $D=D_1\oplus\dots\oplus D_k$ and there is at least one cut
  $D_i$ with $|D_i|\ge |D|$. Now $D_i=D_1\oplus\dots\oplus D_{i-1}\oplus
  D\oplus D_{i+1}\oplus\dots\oplus D_k$ and thus $\mathcal{D}'= \{D\}\cup
  \mathcal{D}\setminus\{D_i\}$ is again a cut basis. By minimality of
  $\mathcal{D}$ we must have $|D_i|=|D|$ and thus $\mathcal{D}'$ is also a
  minimum cut basis. It follows that $D$ is relevant.
\end{proof}

In particular, therefore, all weight-minimal cuts are relevant, see also
\cite{Ren:09}. On the other hand, we can rule out cuts that lead to more
than two connected subgraphs:
\begin{lemma}
  If $D\in\mathfrak{D}$ is relevant, then $D$ is a bond.
  \label{lem:bond}
\end{lemma}
\begin{proof}
  Suppose, for contradiction, that $D$ is not a bond. Then there is a
  cutset $D'\in\mathfrak{D}$ such that $D'\subsetneq D$. Moreover, since
  the cutsets form a vector space, $D''\coloneqq D\oplus D'=D\setminus
  D'\in\mathfrak{D}$. By construction we have $D=D'\oplus D''$ and
  $|D'|<|D|$, and $|D''|<|D|$, and thus $D$ is not a relevant cut.
\end{proof}

The cuts induced by the edges of a Gomory Hu tree are nested and form a
basis of the cut space. The cutset $\partial V'$ of any cut $(V',\bar{V'})$
of $G=(V,E,w)$ therefore can be expressed as a linear combination of the
cut sets defined by a Gomory-Hu tree $T^*$. More precisely,
\begin{equation}
  \partial V' = \bigoplus_{\substack{ uv=e\in F \\ u\in V', v\in\bar{V'}}}
  \partial V_e 
  \label{eq:GHsum}
\end{equation}
where $\partial V_e$ denotes the cutset of the cut $(V_e,\bar V_e)$
obtained by removing $e$ from $T^*$. In other words, $\partial V'$ is the
$\oplus$-sum of the cuts defined by edges in the Gomory-Hu tree $T^*$ with
one vertex in $V'$ and the other one in $\bar{V'}$. Thus $s$ and $t$ are
both in $V'$ or both in $\bar{V'}$ if and only if they are separated by an
even number of cut edges along the unique path $P_{st}$ in $T^*$.

\begin{theorem}
  A cut $\partial V'$ in $G$ is relevant if and only if it is a
  weight-minimal $s,t$-cut for some pair of distinct vertices $s$ and $t$.
  \label{thm:main}
\end{theorem}
\begin{proof}
  \par\noindent\emph{Claim: Every weight-minimal $s,t$-cut is relevant.}
  \par\noindent Let $D^*_{st}$ be an arbitrary but fixed weight-minimal
  $s,t$-cut. Sort the list of all cuts in $G$ by non-decreasing weight such
  that $D^*_{st}$ is the first in the list with its weight. Thus, no
  shorter cut is an $s,t$-cut. Using the canonical greedy algorithm
    \cite{Edmonds:71}, a minimal cut basis $\mathcal{B}$ is obtained by iterating over
    the sorted lists of cuts, adding a cut $D$ to the initial empty set
    $\mathcal{B}$ if and only if $\mathcal{B}\cup\{D\}$ is linearly
    independent. Now consider the set $\mathcal{B}'$, created by the above procedure
    up to the point where all cuts have
    been processed that appear in the list before the first $s,t$-cut
    $D^*_{st}$. By construction, $\mathcal{B}'$ is linearly independent and
    each cutset $\partial U_i\in\mathcal{B}'$ corresponds to a cut
    $(U_i,\bar U_i)$, where we may assume without loosing generality that
    $s,t\in U_i$. By the same argument as in the proof of
    Lemma~\ref{lem:Dst}, $\mathcal{B}'$ cannot generate an $s,t$-cut as an
    $\oplus$-sum.  Thus $\mathcal{B}'\cup\{D^*_{st}\}$ is again linearly
  independent. In the next step, the canonical greedy algorithm
  therefore includes $D^*_{st}$ in the growing minimal weight basis
  $\mathcal{B}'$. It follows that $D^*_{st}$ is
  relevant.\hfill$\triangleleft$
  \smallskip
  
  \par\noindent\emph{Claim: If $\partial V'$ is not a minimum weight
  cut for any pair of vertices $s$ and $t$, then it is not relevant.}

  \smallskip\noindent Let $T^*$ be an arbitrary Gomory-Hu tree for
  $G$. According to Eq. \eqref{eq:GHsum}, $\partial V'$ is a linear
  combination of cutsets $\partial V_e$ defined by $T^*$, where the
    edges $e$ separate vertices $s\in V_e$ and $t\in \bar V_e$.  By (G2)
  and (G3), the weight $\lambda_e$ of the edge $e$ in $T^*$ satisfies
  $\lambda_e=|\partial V_e|$ and, by construction, is the minimum weight of
  cuts separating $s$ and $t$ if $\lambda_e$ is the minimal
    weight \emph{encountered} along the unique path $P_{st}$ in
    $T^*$. Since $\partial V'$ is not a minimum weight cut for any
  pair of vertices, it is in particular not a minimum weight cut for $s$
  and $t$. Thus we have $|\partial V_e|<|\partial V'|$ for every edge $e$
  of the $T^*$ appearing in Eq. \eqref{eq:GHsum}. Thus $D$ is a
  linear combination of strictly shorter cuts. By Lemma~\ref{lem:Vismara},
  therefore, $\partial V'$ it is not relevant.  \hfill$\triangleleft$

  \smallskip\noindent Taken together, a cut $D\in\mathfrak{D}$ is relevant
  if and only if it is a minimum weight cut for at least one pair of
  vertices $s,t\in V$.
\end{proof}

There are graphs with an exponential number of minimal $s,t$-cuts. For
instance the $K_{2,n-2}$ with unit edge-weights and $\{s,t\}$ for one part
of the bipartition \cite{Picard:80}. Since all these cuts are relevant by
Thm.~\ref{thm:main}, the number of relevant cuts, like the number of
relevant cycles \cite{Vismara:97}, may also grow exponentially with the
size of the graph $G$. In the unweighted case, i.e., $w_{pq}=1$ for all
edges $pq\in E$, the number of relevant cuts is bounded by $|E|^{\Delta}$,
where $\Delta$ is the maximum degree because every minimal $s,t$-cut
contains at most $\Delta$ edges and there are no more than $O(|E|^{\Delta})$
distinct edge sets with at most $\Delta$ edges.

This argument does not generalize to arbitrary weights, however: Consider
the graph $G$ obtained from two (potentially identical) instances of binary trees with $n$ leaves by
identifying corresponding leaves. Thus $G$ has $3n-2$ vertices and a
maximum vertex degree of $3$. Suppose all edge weights are unity except for
edges incident with the leaves, which have weight $1/(2n)$. Denoting the
roots of the trees by $s$ and $t$, every cut that contains an edge not
incident with a tree-leaf has weight $>1$, while cutting exactly one of the
two edges incident with a leaf yields an $s,t$-cut with weight $1/2$. There
are $2^n$ distinct cuts of this type, all of which a minimum weight
$s,t$-cut (see Fig. \ref{fig:exp_bintree}). 

\begin{figure}
\centering
\begin{tikzpicture}[yscale=1, every node/.style={draw, circle}, every edge quotes/.style = {above=-4pt, font=\footnotesize, sloped, draw=none}]
\path
(0,0)     node    (1a) {s}
(-2,-1)   node    (2a) {}
+(-1,-1)  node    (3a) {}
+(1,-1)   node (3b) {}
(2,-1)    node    (2b) {}
+(-1,-1)  node (3c) {}
+(1,-1)   node (3d) {}
(3a.center)
+(-.8,-1) node (4a) {}
+(.8,-1)  node (4b) {}
(4a-|1a)  coordinate (M)
(2a|-M)   coordinate (N);

\path
($(4a)+(4b)-(3a)$) node (5a) {}
($(3a)!2!(N)$)     node (6a) {}
($(2a)!2!(N)$)           node (6b) {}
($(1a)!2!(M)$)     node (7a) {t};

\draw 
(1a) edge["1"] (2a) (1a) edge["1"] (2b)
(2a) edge["1"] (3a) (2a) edge["0.1", blue] (3b) 
(2b) edge["0.1", blue] (3c) (2b) edge["0.1", red] (3d)
(3a) edge["0.1", red] (4a) (3a) edge["0.1", red] (4b);
\draw[dashed] 
(5a) edge["0.1", blue] (4a) (5a) edge["0.1", blue] (4b)
(6a) edge["0.1", red] (3b) (6a) edge["0.1", red] (3c)
(6b) edge["1"] (5a) (6b) edge["1"] (6a)
(7a) edge["0.1", blue] (3d) (7a) edge["1"] (6b);
\end{tikzpicture}
\caption{Example of a graph with arbitrary weights and an exponential number of $s,t$-cuts. Following the construction in the main text joining two binary trees (solid vs dashed lines) with $n=5$ leaves, one edge adjacent per leave node has to be included in each minimal cut, yielding $2^5$ possibilities. Two such minimal cuts are marked in blue and red.}
\label{fig:exp_bintree}
\end{figure}

\section{Algorithmic Considerations}

As mentioned in the preliminaries, the set of all minimal $s,t$-cuts,
  which we have shown in Thm.~\ref{thm:main} to coincide with the set
  of relevant cuts, can be computed from a Gomory-Hu tree $T^*$ and the
  directed acyclic graphs introduced by Picard and Queyranne
  \cite{Picard:80} for the edges $uv$ of $T^*$. Both the Gomory-Hu tree
  $T^*$ (for recent advances, see \cite{Abboud:22,Abboud:23}) and the PQ-DAG
  $Q_{s,t}$ for any pair $s,t$ can be computed in polynomial time.

From this combined data structure, all the connectivity cuts, i.e., cuts with capacity equal to the connectivity value,
can be listed with $O(|E| |V|)$ effort \cite{Gusfield:93}. This is feasible as there are at most $\binom{|V|}{2}$ such cuts.
This is reminiscent of the polynomially many ``prototypes'' for relevant
cycles described in \cite{Vismara:97}, from which all relevant cycles can
be listed. There is an important difference, however: For relevant cycles,
there is a polynomial-time algorithm to count their number based on
counting a constrained set of shortest paths \cite{Vismara:97}.  In
contrast, counting the number of minimal $s,t$-cuts is \#P-complete
\cite{Provan:83}, and thus at least as hard as NP-complete decision
problems \cite{Valiant:79}.

By Thm.~\ref{thm:main}, the most direct approach to computing the set of
relevant cuts is to compute, for all distinct pairs $s,t\in V$, the set of
all minimal $s,t$-cuts. We shall see, however, that the
enumeration of all relevant cuts can be sped up considerably with the
help of observations made by Gusfield and Noar in \cite{Gusfield:93}.
First we need some more notation.

%Adapted from \cite{Gusfield:93}.
For any $u, v \in V$, we say that $v$ is a successor of $u$ in $Q_{s,t}$,
if either $u$ and $v$ belong to the same node in $Q_{s,t}$ or if the node
in $Q_{s,t}$ that contains $v$ is a successor of the node in $Q_{s,t}$ that
contains $u$.  Let $Q_{s,t}(X,Y)$ with $X,Y \subseteq V$ and $X \cap Y =
\emptyset$ denote the PQ-DAG resulting from $Q_{s,t}$ by contracting the
set $X$ with all its successors, and contracting the set $Y$ with all its
predecessors. 

\begin{lemma}
Any closed set in $Q_{s,t}(X, Y)$ that contains $s$ but not $t$ is a closed
set in $Q_{s,t}$ that contains $\{s \cup X\}$ but not $\{t \cup Y\}$.
\label{lem:contracted_dag}
\end{lemma}
\begin{proof}
By construction, $s$ is contracted into a super-node with all its successors
and $t$ is contracted into a super-node with all its predecessors to create
$Q_{s,t}$.  Accordingly, the super-nodes containing $s$ and $t$ are the
only vertices with out-degree 0 and in-degree 0, respectively.  For any
vertex $u \in V$, $s$ is a successor of $u$ and $t$ is a predecessors of
$u$, as a consequence of this construction.  Contracting $X$ with all its
successors therefore includes $s$ but not $t$ as a successor. Analogously,
contracting $Y$ with all its predecessors includes $t$ but not $s$ as a
predecessor.  Therefore, the super-node of out-degree 0 in $Q_{s,t}(X, Y)$
represents $\{s \cup X\}$, and the super-node of in-degree 0 represents
$\{t \cup Y\}$.  It follows that any closed set in $Q_{s,t}(X, Y)$ contains
$\{s \cup X\}$ thereby $s$ but not $\{t \cup Y\}$ thereby $t$.
\end{proof}

\begin{lemma} \upshape{\cite[Thm.~2.1]{Gusfield:93}} 
  For a fixed pair of vertices $s,t\in V$, let $(u_1,v_1),\ldots,(u_k,v_k)$
  be the $k$ minimum-weight tree edges along the path between $s$ and $t$
  in $T^*$.  A cut $(V',V'')$ is a minimum $s,t$-cut \textit{if and only if}
  it appears (as a closed set) in at least one PQ-DAG $Q_{u_i,v_i}(s,t)$
  for some $1 \leq i \leq k$, where $Q_{u_i, v_i}(s, t)$ is the graph
  resulting from contracting all successors of $s$ and all predecessors of
  $t$ in $Q_{u_i, v_i}$, respectively, into a single node.
\label{lem:min_path}
\end{lemma}
% \begin{proof}
% See \cite{Gusfield:93} Thm. 2.1.
% \end{proof}

\begin{theorem}
  A cut $D\in\mathfrak{D}$ is \emph{relevant} \textit{if and only if} it is
  represented by at least one closed vertex set in a PQ-DAG $Q_{u, v}$
  derived from an edge in the Gomory-Hu tree $T^*$, hence $\{u,v\} \in
  T^*$.
\end{theorem}
\begin{proof}
``$\Rightarrow$'' By Thm.~\ref{thm:main}, if a cut is relevant, it is a
  weight-minimal $s,t$-cut for some pair of distinct vertices $s$ and $t$.
  By Lemma \ref{lem:min_path} any weight-minimal cut between two vertices
  $s$ and $t$ is represented by a closed set in exactly one PQ-DAG
  $Q_{u,v}(s, t)$, where $uv$ is a weight minimal edge along the path from
  $s$ to $t$ in $T^*$. By Lemma \ref{lem:contracted_dag} any
  $u,v$-cut in $Q_{u,v}(s,t)$ is also a cut in $Q_{u,v}$.

  ``$\Leftarrow$'' Each edge $e=uv$ in a Gomory-Hu tree $T^*$ is
  weight minimal for at least the vertex pair $u$ and $v$ on the path
  consisting of only the edge $uv$ itself. The PQ-DAG $Q_{u,v}$ represents
  all weight minimal cuts between $u$ and $v$ as closed sets. The closed
  sets are therefore weight minimal cuts for two distinct vertices $s$ and
  $t$ and hence relevant as a consequence of Thm.~\ref{thm:main}.
\end{proof}

Therefore, it is sufficient to compute only the PQ-DAGs associated
with edges of a fixed Gomory-Hu tree $T^*$ of $G$.  As
the edges of a Gomory-Hu tree are labeled by the weights of associated
minimal cuts, it is trivial to enumerate the relevant cuts in
non-decreasing order using this method.

In addition to the approach of \cite{Gusfield:93} and its derivation,
several other algorithms for listing relevant cuts are feasible. By
Lemma~\ref{lem:bond}, one could use the bonds of $G$ as a candidate set and
extract relevant cuts using the canonical greedy algorithm similar to
Horton's approach to computing minimum cycle bases \cite{Horton:87}. Very
recently, an algorithm for listing all $(s,t)$-bonds with near linear delay
\cite{Raffaele:24} has become available. A disadvantage is that the set of
bonds must first be generated and then sorted by non-decreasing weight.

A more appealing alternative is to use one of the approaches for ordered
generation of cuts. The first method of this type was proposed in
\cite{Vazirani:92}. Later an improved version based on a different
solutions of maximum flow problems was published \cite{Yeh:10}. Since
  these algorithms produce cuts in the order of non-increasing weights, it
  is not difficult to determine whether a cut $(V',\bar{V'})$ is a minimum
  weight cut for any pair $s\in V$ and $t\in \bar V$. This can be done by
keeping a matrix $M$, initialized by $M_{st}=\infty$ for all $s,t$, and
updated by the following rule: if $M_{u,v}=\infty$, $u\in V$ and $v\in
\bar{V'}$ then $M_{u,v}\leftarrow |\partial V'|$. In particular,
because of Thm.~\ref{thm:main}, it is not necessary to check for
linear independence of a cut with the set of all cuts of strictly
smaller weight as is the case for minimal cycle bases. Moreover, the
generation process can be terminated as soon as (i) minimum weight cuts
have been observed for all vertices, i.e., upon generation of $(V,\bar
V')$, no entries $\infty$ were present in $M$, and (ii) the weight of the
current cut is strictly larger than the weight of the previous one, which
coincides with the largest finite entry in $M$. The enumeration algorithm
in \cite{Yeh:10} has $\tilde{O}(|V| |E|)$ delay, which dominates the
$O(|V|^2)$ effort for ``filtering out'' cuts that are not minimal
$u,v$-cuts for any pair $u,v$.

\section{Computational Experiments} 

\begin{figure}
\centering
\includegraphics[width=0.9\textwidth]{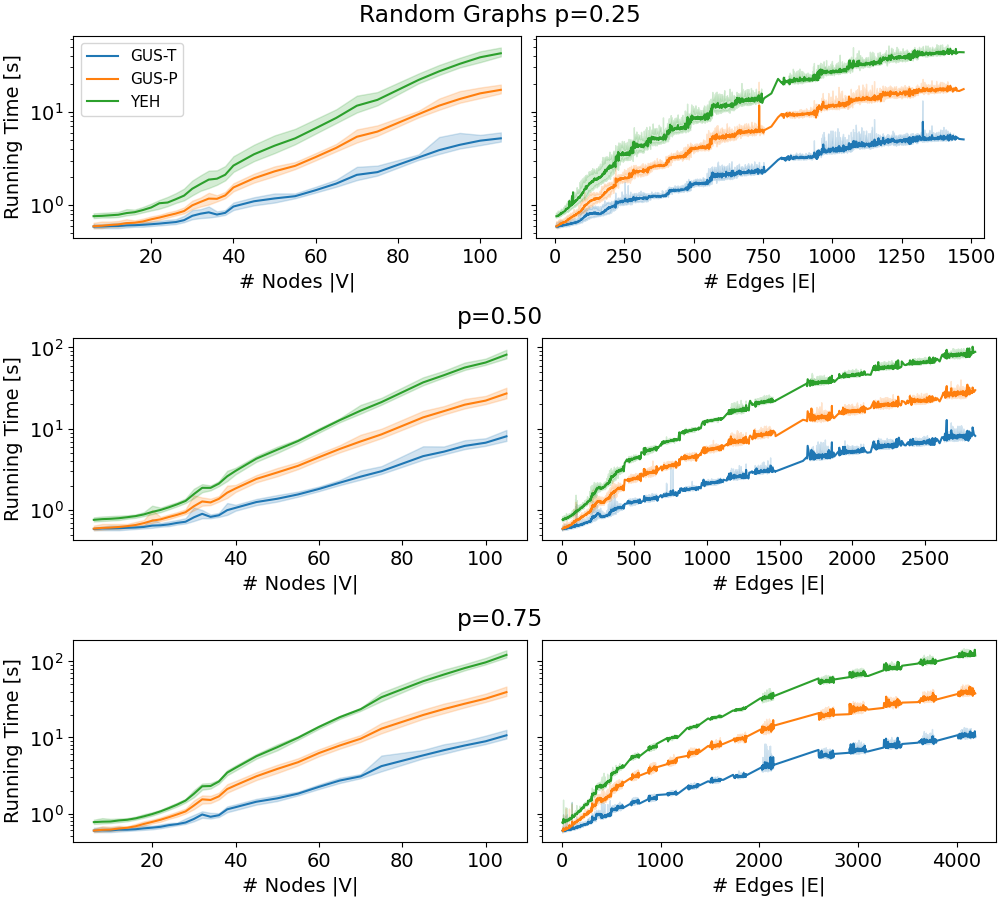}
\caption{Algorithmic performance for random graphs, partitioned by edge
  probability in the Erdős–Rényi model. Running times were tracked
  and related to the number of nodes (left) and the number of edges
  (right) of graphs. Lines represent the average running time, while
  environments show the 95\% percentile interval. All plots present
  logarithmic y-scales. Regardless of edge density, methods show similar
  pattern.}
\label{fig:plot_random}
\end{figure}

\begin{figure}
\centering
\includegraphics[width=0.9\textwidth]{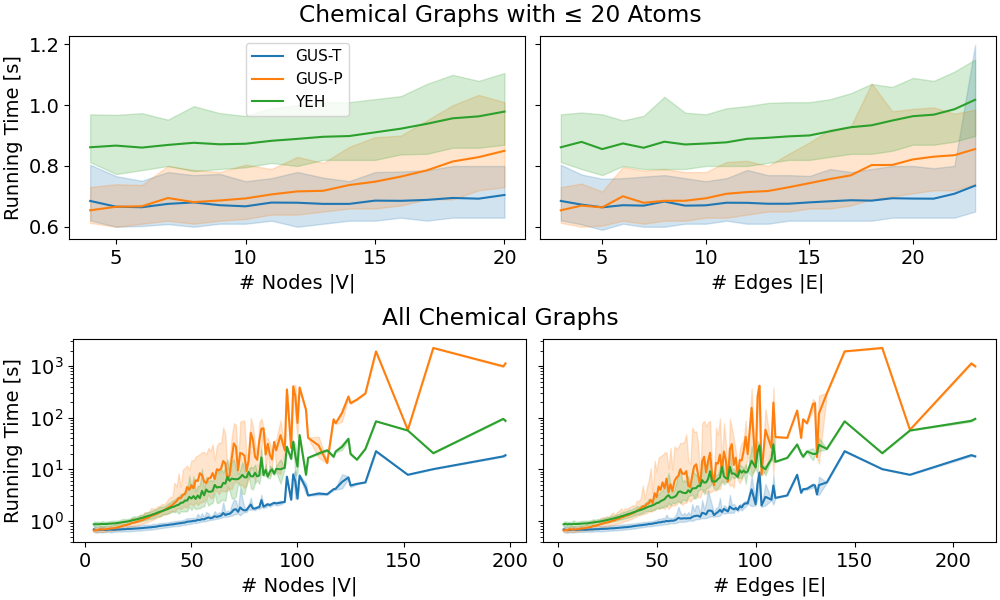}
\caption{Algorithmic performance for chemical graphs. Running times
  were tracked and related to the number of nodes (left) and the
  number of edges (right) of graphs. For improved resolution, chemical
  graphs containing 20 or less atoms are plotted independently (above),
  contrasting the full dataset (below). Lines represent the average
  running time, while environments show the 95\% percentile
  interval. Please note that only the lower graphs follow a logarithmic
  y-scale.}
\label{fig:plot_chem}
\end{figure}

\begin{table}[bht]
\caption{Summary of benchmarked chemical graphs. \texttt{GUS-T} performs
  best for the vast majority of cases, with exception of few very small
  graphs.}  
\begin{center}
\begin{tabular}{||c || c c c||} 
 \hline
 \textbf{\textit{Category}} & \texttt{GUS-T} & \texttt{GUS-P} & \texttt{YEH} \\ [0.5ex] 
 \hline\hline
Total Running Time [s]     & 8340.50 & 35035.97 & 18616.66 \\
Fastest Method [\# Graphs] & 9916    & 84       & 0        \\
 \hline
\end{tabular}
\label{tbl:bench} 
\end{center}
\end{table}

We implemented both methods derived from \cite{Gusfield:93}, enumerating
closed sets in PQ-DAGs for all pairs of vertices (\texttt{GUS-P}) and only
for PQ-DAGs for edges of $T^*$ (\texttt{GUS-T}), as well as the algorithm
for enumerating cuts by non-increasing weights proposed by \cite{Yeh:10}
(\texttt{YEH}).  All methods where implemented in \textit{Python}
  (v3.11) utilizing the ubiquitous graph library \textit{NetworkX}
  (v. 3.3).  All benchmarks where run on a consumer PC (Intel i9-9900K)
  with sufficient memory and without parallelization.  Memory usage overall
  was negligible for all tested instances and is therefore not reported.

Benchmarks were performed on both artificial and realistic data.  As
generating cuts is of particular interest in the domain of computational
chemistry, we used $10,000$ individual chemical graphs from the USPTO-10K dataset \cite{morin2023molgrapher},
a curated subset of reactions reported in the United States Patent and Trademark Office.
Graphs ranged in size between 4 and 198 nodes with a median of 30 at 3 to 211
edges and a median of 33. As such, chemical graphs are considered to be
sparse.  Random graphs were constructed following the classic Erdős–Rényi
model \cite{Erdos:1959:pmd}, generating 1000 valid instances per parameter variant, i.e. per each
6 to 105 nodes (at step size 2 until 40, then at step size 5) and per edge
probability 0.25, 0.5, and 0.75. Only connected graphs where considered.
The cut basis was computed for each method and dataset, logging individual
running times.  Times were averaged over $3$ iterations for chemical
graphs.

Random graphs and chemical graphs both show a consistent, but markedly
different hierarchy between methods.  On all random graphs (see
Fig.~\ref{fig:plot_random}), \texttt{GUS-T} consistently outperformed both
\texttt{GUS-P} and \texttt{YEH}, while \texttt{GUS-P} consistent
outperformed \texttt{YEH}. Margins widen at increasing edge and node
counts.  Edge probability does not seems to have any influence on this
relationship.
 
In chemical graphs, unsurprisingly, \texttt{GUS-T} also consistently
outperformed both \texttt{GUS-P} and \texttt{YEH} for all but very small
graphs (see Fig. \ref{fig:plot_chem}, Tbl.~\ref{tbl:bench}). However,
consistently, but more pronounced for larger graphs, \texttt{YEH} was
significantly faster than \texttt{GUS-P}.  Although \texttt{GUS-P} was the
fastest algorithm for almost 1\% of the graphs, manual inspection showed
that this observation is most likely an artifact caused by fluctuations in
loading times for very small graphs.  Overall, \texttt{YEH} was a factor
2.2 slower then \texttt{GUS-T} to create a cut basis. 

\section{Concluding Remarks}

With \texttt{GUS-T}, we derived an effective algorithm to compute exactly all relevant cuts of a graph, based on two central observations: relevant cuts coincide with the cuts that are weight-minimal for some pair of vertices, which in turn are fully encoded as the cuts implied by the edges of a Gomory Hu tree. Picard-Queyranne DAGs yield an effective framework towards this goal. \texttt{GUS-T} consistently outperformed other methods in our benchmarks, regardless of graph properties. Nevertheless, it should be acknowledged that \texttt{YEH} is overall more
versatile when enumerating cuts.

\texttt{GUS-T} is in particular well suited to enumerate relevant cuts in chemical applications, where cuts relate to chemical reaction mechanism, such as for cut-and-join crossover operators \cite{Domschke:25a}.
It remains an open question, however, whether the restriction of cuts to the set
of relevant cuts remains sufficient to ensure global reachability in the space
of chemical graphs in conjunction with degree-preserving or even more
restricted join operations.

\section*{Availability}

A reference implementation of the methods developed in this paper are
available via GitHub \url{https://github.com/TGatter/Cut-Base-Toolkit}.

\section*{Author Contributions}

ND, TG, RG, and PFS designed the study, obtained the theoretical results and
drafted the manuscript. ND and TG created the reference implementation and
executed the benchmarks.  The students of the Advanced Bioinformatics
Course 2024 implemented \cite{Yeh:10} as a part of a course project.

\subsection*{Acknowledgements} 
This work was supported by the German Research Foundation (DFG) within
SPP2363 (460865652) and the Novo Nordisk Foundation as part of MATOMIC
(0066551). Research in the Stadler lab is supported by the German Federal
Ministry of Education and Research BMBF through DAAD project 57616814
(SECAI, School of Embedded Composite AI). We gratefully acknowledge the
contributions of the students of the Advanced Bioinformatics Course 2024:
Marie Freudenberg, Julius Guntrum, Florian K{\"u}hn, Lukas Marche, Louisa
Marie von Menges, and Kiet Ngo Tuan.

\bibliographystyle{adamjoucc}
\bibliography{relcut}
\label{sect:biblio}

\end{document}